\documentclass[12pt,twoside]{amsart}

\usepackage{amssymb}

\newtheorem{thm}{Theorem}

\newtheorem{lem}[thm]{Lemma}
\newtheorem{cor}[thm]{Corollary}

\newtheorem{conj}[thm]{Conjecture}
   
\theoremstyle{definition}

\newtheorem{say}[thm]{}

\newtheorem{ack}{Acknowledgments}

\newtheorem{defn-thm}[thm]{Definition--Theorem}  

\theoremstyle{remark}

\renewcommand{\c}[0]{{\mathbb C}}  

\renewcommand{\o}[0]{{\mathcal O}} 

\newcommand{\p}[0]{{\mathbb P}}

\newcommand{\q}[0]{{\mathbb Q}}
\newcommand{\map}[0]{\dasharrow}
\newcommand{\qtq}[1]{\quad\mbox{#1}\quad}

\newcommand{\mult}[0]{\operatorname{mult}}

\newcommand{\Hom}[0]{\operatorname{Hom}}

\newcommand{\aut}[0]{\operatorname{Aut}}

\newcommand{\chow}[0]{\operatorname{Chow}}

\begin{document}
\bibliographystyle{amsplain}

\title{Characterizations of $\p^n$ in arbitrary characteristic}
\author{Yasuyuki Kachi  and J\'anos Koll\'ar}

\maketitle

There are several results which characterize $\p^n$
among projective varieties. The first substantial theorem
in this direction is due to 
\cite{HiKo}, using the topological properties of $\c\p^n$.
Later characterizations approached 
 $\p^n$
as the projective variety with the most negative canonical bundle.
The first such theorem is due to \cite{kob-o} and later it was generalized
in the papers  \cite{ion, fuj}. A  characterization using the tangent  bundle
was given in \cite{mori, s-y}.
Other characterizations using vector bundles
are given in the papers \cite{y-z, sato, kachi}.

The results of \cite{kob-o, ion, fuj} 
 ultimately relied on the Kodaira vanishing
theorem, thus they were restricted to characteristic zero.
The aim of this note is to provide an argument  which does not use higher
cohomologies of line bundles. This makes the proofs work
in arbitrary characteristic. Thus
(\ref{main.thm}) generalizes \cite{ion, fuj}
and (\ref{main.cor}) is a characteristic free version of  \cite{kob-o}. 
The proofs of \cite{ion, fuj} work if we know enough about  extremal
contractions of smooth varieties. Thus the results have been known
in any characteristic for surfaces and they can be read off
from the classification of contractions of smooth 3--folds
in positive characteristic
developed in the papers \cite{koll-pcone, shep-b, megy}.

For the rest of the paper  all varieties are
over  an algebraically closed field
of arbitrary characteristic.

\begin{thm}\label{main.thm}   Let $X$ be a smooth projective variety 
of dimension $n$  and $H$ an ample divisor on $X$.
 Then the pair $(X,H)$
satisfies one of the following:
\begin{enumerate}
\item $(n-1)H+K_X$ is nef.
\item $X\cong \p^n$ and $H$ is a hyperplane section.
\item $X$ is a quadric in $ \p^{n+1}$ and $H$ is a hyperplane section.
\item $X$ is the projectivization of a rank $n$  vector bundle over
 a smooth curve $A$ and $H=\o(1)$.
\end{enumerate}
\end{thm}

With  a little bit of work on the scrolls in (\ref{main.thm}.4)
this implies the following
 characterization  of projective spaces and hyperquadrics.
 
\begin{cor}\label{main.cor}
 Let $X$ be a smooth projective variety of dimension $n$.
\begin{enumerate}
\item  $X$ is isomorphic to $\p^n$ iff there is an ample divisor $H$ such that
$-K_X$ is numerically equivalent to $(n+1)H$.
\item  $X$ is isomorphic to a hyperquadric $\q^n\subset \p^{n+1}$ 
iff there is an ample divisor $H$ such that
$-K_X$ is numerically equivalent to $nH$.\qed
\end{enumerate}
\end{cor}

Our proof of (\ref{main.thm}) also yields the following:

\begin{cor}\label{2nd.cor}
 Let $X$ be a smooth projective variety of dimension $n$.
Assume that  $K_X$ has negative intersection number with some curve. 
Then  $X$ is isomorphic to $\p^n$ iff 
\begin{enumerate}
\item $(-K_X\cdot C)\geq n+1$ for every rational curve $C\subset X$, and
\item  $(-K_X)^n\geq (n+1)^n$.
\end{enumerate}
\end{cor}

It is conjectured that (\ref{2nd.cor}.1) alone characterizes $\p^n$.
A positive solution in characteristic zero is announced in
\cite{c-m}. It has also been conjectured that
(\ref{2nd.cor}.2) also characterizes $\p^n$ among Fano varieties.
This, however, turned out to be false  for $n\geq 4$ \cite{bat}. 
\medskip

 The first step  of the proof  of (\ref{main.thm})
is to identify the family of
lines on $X$. The results of \cite{mori} show that there
is a family of  rational curves  $\{C_t\}$ on $X$
such that $(H\cdot C_t)=1$ and this family has the expected dimension.
Focusing on these  rational curves  leads to 
other characterizations of $\p^n$ and of $\q^n$, which
 strengthen (\ref{main.cor})  in characteristic zero.

\begin{thm}\label{abw-thm}
 Let $X$ be a smooth projective variety 
of dimension $n$
over $\c$ and
  $L$  an ample line bundle  on $X$.
\begin{enumerate}
\item \cite{abw} Assume that for every
$x,y\in X$ there is a rational curve $C_{xy}$ through $x,y$ 
with $(L\cdot C_{xy})=1$. Then  $X\cong \p^n$ and $L\cong \o(1)$.
\item \cite{ks} Assume that for every
$x,y,z\in X$ there is a rational curve $C_{xyz}$ through $x,y,z$ 
with $(L\cdot C_{xyz})=2$. Then  $X\cong \p^n$ 
or $X\cong \q^n$ and $L\cong \o(1)$.
\end{enumerate}
\end{thm}

The positive characteristic version of  (\ref{abw-thm})
is still open.
Concentrating on higher degree  curves
leads to the following characterization of Veronese varieties,
which gives another generalization of (\ref{main.cor}.1)
and also leads to a weaker version of (\ref{main.thm}).

\begin{thm}\cite{kachi}  Let $X$ be a smooth projective variety 
of dimension $n$ (over an algebraically closed field of arbitrary 
characteristic)
 and
  $L$  an ample divisor  on $X$. Assume that:
\begin{enumerate}
\item For general $x,y\in X$ there is a connected (possibly reducible) curve
$C_{xy}$ through $x,y$ 
with $(L\cdot C_{xy})\leq d$.
\item For every irreducible subvariety
$Z\subset X$ of dimension at least 2  we have
$(L^{\dim Z}\cdot Z)\geq d^{\dim Z}$. 
\item Fujita's sectional genus of $\frac1{d}L$ is less than 1, that is,
$$
(L^{n-1}\cdot (\textstyle{\frac{n-1}{d}}L+K_X))< 0.
$$
\end{enumerate}
Then $X\cong \p^n$ and $L\cong \o(d)$.
\end{thm}

\begin{say}[Proof of (\ref{main.thm})]{\ }
As we noted, we may  assume that $n\geq 3$.
 We are done if  $K_X+(n-1)H$ is  nef.  

If  $K_X+(n-1)H$ is not  nef, then
by the cone theorem of Mori (cf. \cite[III.1]{ratbook}) there is a $c>0$
and a rational curve $C\subset X$ 
such that
\begin{enumerate}
\item $K_X+(n-1+c)H$ is nef,
\item $C\cdot (K_X+(n-1+c)H)=0$, and
\item $-(C\cdot K_X)\leq n+1$.
\end{enumerate}
Thus we see that 
$$
(C\cdot H)\leq \frac{-(C\cdot K_X)}{n-1+c}\leq \frac{n+1}{n-1+c}<2.
$$
Hence  $C\cdot H=1$ and $-(C\cdot K_X)=n$ or $n+1$.
Let $f:\p^1\to X$ be the normalization of $C$, $0\in \p^1$ a 
point and  $x=f(0)$. 
Then by \cite[II.1]{ratbook}
$$
\begin{array}{ll}
\dim_{[f]}\Hom(\p^1,X)& \geq -(C\cdot K_X)+n,\qtq{and}\\
\dim_{[f]}\Hom(\p^1,X, 0\mapsto x)& \geq -(C\cdot K_X).
\end{array}
$$

We consider 3 separate cases corresponding to the three outcomes
 (\ref{main.thm}.2-4). 
\end{say}

\begin{say}[Case 1: $-(C\cdot K_X)=n+1$]{\ }\label{case1}

In order to go from homomorphisms to curves in $X$, we have to
quotient out by the automorphism group. Thus 
$\Hom(\p^1,X, 0\mapsto x)$ gives an $(n-1)$-dimensional
family of rational curves through $x\in X$
since $\aut(\p^1,0)$ is 2--dimensional.
This implies that the Picard number of $X$ is 1 
(cf.\ \cite[IV.3.13.3]{ratbook}). Thus $-K_X$ is numerically equivalent 
to $(n+1)H$. 
The following lemma shows
that $X\cong \p^n$. We formulate it for singular varieties,
the setting  needed for the third case.
\end{say}

\begin{lem}\label{p^nlem}
 Let $X$ be a normal, projective variety and $x\in X$
a smooth point. Let $H$ be an ample  Cartier divisor on $X$
and let $\{C_t:t\in T\}$
 be an $(n-1)$-dimensional family of curves through $x$.
Assume that $(C_t\cdot H)=1$ and $(-K_X\cdot H^{n-1})>(n-1)(H^n)$.
Then $X\cong \p^n$, $H=\o(1)$ and the $C_t$ are lines through $x$.
\end{lem}

Proof. By Riemann--Roch,
$$
h^0(X, \o_X(mH))=m^n\frac{(H^n)}{n!}+
   m^{n-1}\frac{-(K_X\cdot H^{n-1})}{2(n-1)!}+O(m^{n-2}).
$$
For a section of a line bundle it is
$$
\binom{m+n-1}{n}=m^n\frac{1}{n!}+
   m^{n-1}\frac{n-1}{2(n-1)!}+O(m^{n-2})
$$
conditions to vanish at  a smooth point $x\in X$ to order $m$.
Comparing the two estimates we see that for large $m$
$$
\begin{array}{lll}
H^0(X, \o_X(mH)\otimes I_x^{m+1})&\geq
(\mbox{const})\cdot m^n &\qtq{if $(H^n)>1$, and}\\
H^0(X, \o_X(mH)\otimes I_x^m)&\geq
(\mbox{const})\cdot m^{n-1} &\qtq{if $(H^n)=1$.}
\end{array}
$$
Pick any member $D\in |\o_X(mH)\otimes I_x^{m+1}|$
(resp.\ $D\in |\o_X(mH)\otimes I_x^m|$)
and take $C_t\not\subset D$. Then
$(C_t\cdot D)=m$. On the other hand,
$(C_t\cdot D)\geq \sum_{y\in C_t}\mult_yD\geq m+1$ (resp.\ $\geq m$). 
This is a contradiction if $(H^n)>1$. If $(H^n)=1$ then 
 $(C_t\cap D)=\{x\}$ and the linear system
$|\o_X(mH)\otimes I_x^m|$ is constant along $C_t$.

By varying $x\in C_t$, we also see that   the line bundle
$\o_{C_t}(mH)$ has a section with an $m$-fold zero
at a general point of $C_t$. Since $\deg \o_{C_t}(mH)=m$,
this implies that a general $C_t$ is a smooth rational curve.

Set $Y=B_xX$ with exceptional divisor $E$ and projection
$\pi:Y\to X$. Set $M=\pi^*H-E$. Then $H^0(\o_Y(mM))\cong
H^0(X,\o_X(mH)\otimes I_x^m)$, hence 
$h^0(\o_Y(mM))\geq (\mbox{const})\cdot m^{n-1}$.
For large $m$, let $p:Y\map Z$ denote the map induced by
$|mM|$. Up to birational equaivalence, $p$ does not depend on
$m$ and  it has connected 1--dimensional general fibers.
Moreover, 
 for any subset $V\subset Y$ of dimension at most $n-2$
there is a member $D\in |mM|$ containing $V$.

Every curve $C_t$ lifts to a curve $C'_t$ on $Y$
and there is an open subset $T^0\subset T$ such that the curves
$\{C'_t:t\in T^0\}$ form a single algebraic family.
(In our case these correspond to the curves 
$C_t$ 
which are smooth at $x$.)
By shrinking $T$ we pretend that $T=T^0$.
Note that $(C'_t\cdot E)=1$ and $(C'_t\cdot M)=0$.

We may assume that $Z$ is proper. Let $V\subset Y$ be the set where $p$ 
is not defined and choose
$V\subset D\in |mM|$. 

 If a curve   $C'_t$ passes through $V$ 
then $C'_t\subset D$ since $(M\cdot C'_t)=0$. The
curves $C'_t$ cover an open subset of $Y$, hence  a general $C'_t$
is disjoint from $D$. Thus $p$ is defined everywhere along a general $C'_t$
and $C'_t$
is a fiber of $p$ (at least set theoretically). 
 Hence there
 are open subsets $Y^0\subset Y$  and $Z^0\subset Z$
such that $p^0:Y^0\to Z^0$ is proper and flat.
$(E\cdot C'_t)=1$, thus $E$ is a rational section
of $p$. In particular, a general fiber of $p$ is reduced. 
Thus, possibly after shrinking $Y^0$, we may assume that  $p^0$ is a
 $\p^1$-bundle.  Set $E^0:=E\cap Y^0$.

$Z^0$ parametrizes curves in $Y$ and by looking at the
image of these curves in $X$
we obtain a morphism $h:Z^0\to \chow(X)$. 
(See \cite[I.3--4]{ratbook} for the definition of $\chow$
and its basic properties.)
Let $Z'$ be the
normalization of the closure of the image of $Z^0$ in $\chow(X)$, $p':U'\to Z'$
the universal family and $u':U'\to X$ the cycle map.
We want to prove that $U'=Y$.

$h$ induces a morphism between the universal families
$h':Y^0\to U'$ which sits in a diagram
$$
Y^0\stackrel{h'}{\to} U'\stackrel{u'}{\to} X\stackrel{\pi^{-1}}{\map} Y.
$$
The composite is an open immersion hence birational. 
Thus the intermediate maps are
also birational. $h'$ is thus an isomorphism near a general
fiber of $p^0$, and so $p':U'\to Z$ is a $\p^1$-bundle over an open set.
$Z'$ parametrizes curves wich have intersection number 1 with an ample
 divisor. Thus $Z'$ parametrizes irreducible curves with multiplicty 1.

Let $E'\subset U'$ be the closure of $h'(E^0)$. 
$E'$ is a rational section of $p'$ and $u'(E')=\{x\}$. 
If $E'$ contains a whole fiber of $p'$
then  $u'(E')=\{x\}$ implies that every fiber is mapped to a point
which is impossible. So,
 $E'$ is a section of $p'$.

Let $B\subset U'\setminus E'$  be a curve such that $u(B)\in X$
is a point. Then $S:=(p')^{-1}(p'(B))$ is a ruled surface
and $u'$ contracts the curves  $B$ and $E'\cap S$
to points. By bend--and--break (cf. \cite[II.5.5.2]{ratbook})
this leads to a curve  $C^*\subset X$
such that $(C^*\cdot H)<1$, but this is impossible. Thus
$u'$ is quasifinite on $U'\setminus E'$.
$u'$ is also birational, hence
$$
u':U'\setminus E'\to X \qtq{is an open embedding.}
$$
Its image contains $X\setminus \{x\}$ and it is not projective, hence 
$U'\setminus E'\cong X\setminus \{x\}$.

Let us now look at the birational map $\phi:={u'}^{-1}\circ \pi:Y\map U'$. 
Both $U'$ and $Y$ contain $Y^0$ and $X\setminus\{x\}$ as  open sets, and
$\phi$ is the identity on $Y^0$ and $X\setminus\{x\}$.
 Thus $\phi$ is
an isomorphism outside the codimension 2 sets
$E\setminus E^0$ and $E'\setminus E^0$. 
Since   $Y/X$ has relative Picard number 1,
this implies that $\phi$ is a morphism 
by (\ref{mats-mumf}). 
Thus $\phi$ is  an  isomorphism since $\phi$ can not contract
a subset of $E\cong \p^{n-1}$ without contracting $E$.  

In particular $E'\cong\p^{n-1}$.
$p':U'\to E'$ is flat with reduced fibers
by \cite[Ex.III.10.9]{hartsh}, hence it is a $\p^1$-bundle.
Thus $X\cong \p^n$ by an
easy argument (see, for instance,  \cite[V.3.7.8]{ratbook}).\qed
\medskip

The following lemma is essentially in 
\cite{mat-mum}.

\begin{lem}\label{mats-mumf} Let $Z_i\to S$ be projective morphisms
with $Z_1$ smooth.
Let $\phi:Z_1\map Z_2$ be a birational map 
of $S$-schemes
which is an isomorphism outside
the codimension 2 subsets $E_i\subset Z_i$. Assume that
the relative Picard number $\rho(Z_1/S)$ is 1
and $Z_1\setminus E_1\to S$ is not quasifinite.
 Then $\phi$ is a morphism.
\end{lem}

Proof. We may assume that $S$ is affine.
Let $H_2$ be a relatively ample divisor on $Z_2$ and $H_1$ its 
birational transform. $H_1$ is also relatively ample
because $\rho(Z_1/S)=1$ and $-H_1$ can not be relatively nef
(since $H_1$  is effective when restricted to a positive dimensional fiber of
$Z_1\setminus E_1\to S$). 
Thus $|mH_2|$ and $|mH_1|$ are both base point free for $m\gg 1$
and these are the birational transforms of each other by $\phi$.
Let $\Gamma\subset Z_1\times_S Z_2$ be the closure of the graph of $\phi$
with projections  $\pi_i: \Gamma\to Z_i$. If
$\pi_1^{-1}(z)$ is positive dimensional then every member of
$|mH_2|$ intersects $\pi_2(\pi_1^{-1}(z))$, thus
$|mH_1|=\phi^{-1}|mH_2|$ has $z$ as its base point, a contradiction.
Thus $\pi_1$ is an isomorphism and $\phi$ is a morphism.\qed

\begin{say}[Proof of (\ref{2nd.cor})]{\ }

By \cite{mori} there is  a rational curve $C\subset X$ such that
$0<(-K_X\cdot C)\leq n+1$. 
Fix an ample divisor $L$ and pick a rational curve $C\subset X$ such that
$0<(-K_X\cdot C)$ and $(L\cdot C)$ is minimal. 
Then  $(-K_X\cdot C)= n+1$
by (\ref{2nd.cor}.1) and $C$ is not numerically equivalent
 to a reducible curve  whose irreducible components are rational.
Thus, as in (\ref{case1}),  we obtain that
the Picard number of $X$ is 1 and $-K_X$ is ample.
Set $H:=-\frac{1}{n+1}K_X$. $H$ is a $\q$-divisor
such that  $(H^n)\geq 1$ and $(H\cdot C)\geq 1$
for every rational curve $C\subset X$. It is easy to see
that the proof of (\ref{p^nlem}) works for $\q$-divisors
 satisfying these
assumptions. Thus we obtain (\ref{2nd.cor}).\qed
\end{say}

Returning to the proof of (\ref{main.thm}),
we are left with the cases when $-(C\cdot K_X)=n$. 
Then we have an at least $(n-2)$-dimensional
family of curves through every $x\in X$. 
If there is an $x\in X$ such that all curves of degree 1 through $x$ cover 
$X$,  then we
obtain that $X\cong\p^n$ and this leads to a contradiction.
Thus
 for every $x\in X$ the  curves $\{C_t:x\in C_t\}$
  sweep out a divisor
 $B_x$. The divisors form an algebraic family for 
$x$ in a suitable open set $X^0\subset X$.

\begin{say}[Case 2: $-(C\cdot K_X)=n$ and
 $B_{x_1}\cap B_{x_2}\neq \emptyset$ for  $x_1,x_2\in X^0$]{\ }

By assumption, any two points of $X$ are connected
by a chain of length 2 of the form $C_{t_1}\cup C_{t_2}$,
thus 
 $ X$ has Picard number 1 by \cite[IV.3.13.3]{ratbook}.
 Thus we see that
$-K_X\equiv nH$. The computation of the genus of a general 
complete intersection curve of members of $mH$ for odd  $m\gg 1 $ 
shows that $(H^n)$ is even. Thus $(H^n)\geq 2$.

Pick  general $x_1,x_2$ 
such that there is no degree 1 curve through $x_1$ and $x_2$ and pick
 curves $C_i\supset x_i$ such that 
$C_1\cap C_2\neq \emptyset$. We can view $C_1\cup C_2$
as the image of a map
$f:\p^1\vee \p^1\to X$ where $\p^1\vee \p^1$ denotes
the union of 2 lines in $\p^2$.
Let $y_i\in \p^1\vee \p^1$ be a preimage of $x_i$. By
the usual estimates (cf.\ \cite[II.1.7.2]{ratbook}) we obtain that
$$
\dim_{[f]}\Hom(\p^1\vee \p^1,X, y_i\mapsto x_i)\geq 2n+n-2n=n.
$$
On the other hand, these maps correspond to pairs
of curves $C'_1\cup C'_2$ where $C'_i$ passes through $x_i$
and $C_1\cap C_2\neq \emptyset$. By our assumption
these form an $n-2$-dimensional family. The automorphism group
$\aut(\p^1\vee \p^1, y_1,y_2)$ accounts for the missing 2 dimensions.
Viewing $\p^1\vee \p^1$ as a reducible plane conic,
\cite[II.1.7.3]{ratbook} implies that there is an at least
$(n-1)$-dimensional family of degree 2 rational curves
  $\{A_s\subset X:s\in S\}$ 
which pass through both of $x_1, x_2$. 

As in Case 1,  we obtain that
$$
\begin{array}{lll}
H^0(X, \o_X(mH)\otimes I_{x_1}^{m+1}\otimes I_{x_2}^{m+1})&\geq
(\mbox{const})\cdot m^n &\qtq{if $(H^n)>2$, and}\\
H^0(X, \o_X(mH)\otimes I_{x_1}^m\otimes I_{x_2}^m)&\geq
(\mbox{const})\cdot m^{n-1} &\qtq{if $(H^n)=2$.}
\end{array}
$$
Pick any member $D\in |\o_X(mH)\otimes I_{x_1}^{m+1}\otimes I_{x_2}^{m+1}|$
(resp. $D\in |\o_X(mH)\otimes I_{x_1}^m\otimes I_{x_2}^m|$)
and take $A_s\not\subset D$. Then
$(A_s\cdot D)=2m$. On the other hand,
$(A_s\cdot D)\geq \sum_{y\in A_s}\mult_yD\geq 2m+2$ (resp.\ $ \geq 2m$). 
Thus
$(H^n)=2$,  $(A_s\cap D)=\{x_1,x_2\}$ and the linear system
$|\o_X(mH)\otimes I_{x_1}^m\otimes I_{x_2}^m|$ is constant along $A_s$
for general $s$. 

We also see that  for general $A_s$, the line bundle
$\o_{A_s}(mH)$ has a section with an  $m$-fold zero
at two general points of $A_s$. Since $\deg \o_{A_s}(mH)=2m$,
this implies that a general $A_s$ is a smooth rational curve.

Set $Y=B_{x_1x_2}X$ with exceptional divisors $E_1,E_2$ and projection
$\pi:Y\to X$. Set $M=\pi^*H-E_1-E_2$. Then $|mM|=
|\o_X(mH)\otimes I_{x_1}^m\otimes I_{x_2}^m|$.
As in Case 1, we obtain an open set $Y^0\subset Y$ and
a $\p^1$-bundle  $Y^0\to Z^0$ with two sections
$E_i\cap Y^0$. 
Construct $p':U'\to Z'$ and $u:U'\to X$
using $\chow(X)$ as before.  
Let $z\in Z'$ be any point and  $A_z\subset X$ the 1--cycle
 corresponding to $z$. $A_z$ has degree 2 and it
 passes through the points $x_1,x_2$. We have assumed that 
there is no degree 1 curve passing through  $x_1,x_2$,
thus either $A_z$ is irreducible and $(H\cdot A_z)=2$
or $A_z$ has 2 different irreducible components $A_z^1\cup A_z^2$
and $(H\cdot A_z^i)=1$.

 The  key point is to establish that
the cycle map $u$ is quasifinite on
$U'\setminus (E'_1\cup E'_2)$.
Assuming the contrary, we have a normal surface $S$ 
and  morphisms $p':S\to Z^*\subset Z'$ and $u:S\to X$
such that every fiber of $p':S\to Z^*$ is either $\p^1$
or the union of 2 copies of $\p^1$ and $u$
contracts $E_i\cap S$ and another curve $B$ to points.

 A version
of bend--and--break (\ref{3ptb&b}) establishes that in this case
$u({p'}^{-1}(z))$ is independent of $z\in Z^*$.
This is, however, impossible since $Z'\to \chow(X)$
is finite.

Again using (\ref{mats-mumf}) we conclude that
$U'\cong Y$ and $Z'\cong E'_i\cong \p^{n-1}$. 
Therefore
$$
\o_Y(\pi^*H-E_1-E_2)\cong {p'}^*\o_{E'_i}(1)\cong {p'}^*\o_{\p^{n-1}}(1)
$$
is generated by global sections.
Pushing these sections down to $X$ and varying the points $x_i$
we obtain that $\o_X(H)$ is generated by $n+2$ global sections.
Thus $|H|$ gives  a finite morphism $X\to \p^{n+1}$
whose image is a quadric since $(H^n)=2$.  Thus $X$ is a smooth quadric.

\end{say}

\begin{lem} [3 point bend--and--break]\label{3ptb&b}
  Let  $S$ be a normal and proper
 surface which is a conic bundle (that is, there is a morphism
$p:S\to A$  such that every (scheme theoretic) fiber is
isomorphic to a plane conic).
 Let $E_1,E_2\subset S$ be disjoint sections of $p$
and $B\subset S$ a multisection. Assume that
 every singular fiber $F_a$ of $p$ has 2 components $F_a=F_a^1\cup F_a^2$,
and $(F_a^i\cdot E_i)=0$ for $i=1,2$. 
Let $L$ be a nef divisor on
 $S$ such that $(L\cdot F_a^i)=1$ for every $F_a$ and $i$ and
 $(L\cdot E_1)=(L\cdot E_2)=(L\cdot B)=0$. Then $S\cong A\times \p^1$,
the $E_i$ are flat sections and $B$ is a union of flat sections.
\end{lem}

Proof. The  Picard group of  $S$ is generated by the classes
$F_a^i, E_1$ with rational coefficients (cf. \cite[IV.3.13.3]{ratbook}).
  From the conditions
$(L\cdot F_a^i)=1$ we conclude that $L=E_1+E_2+dF$ for some $d$
(and  $(E_1^2)=(E_2^2)$). 
$(L^2)\geq 0$ and  $(L\cdot E_i)=0$, so by 
the Hodge index theorem $(E_i^2)\leq 0$. Hence we conclude that
$d=-(E_i^2)\geq 0$. This implies that
$$
(L\cdot B)=(E_1\cdot B)+(E_2\cdot  B)+\deg (B/A)\cdot (-E_1^2),
$$
All terms on the right hand side are nonnegative.
Thus $(E_1^2)=(E_2^2)=0$. Since  $(E_1\cdot E_2)=0$, this implies that
$E_1$ and $E_2$ are algebraically equaivalent.
 Thus there are no singular fibers.
\qed

\begin{say}[Case 3: $-(C\cdot K_X)=n$ and   $B_{x_1}\cap
 B_{x_2}=\emptyset$ for general  $x_1,x_2$]

Take a general $x_1\in X$ and a general $x_3\in B_{x_1}$.
If $\dim (B_{x_1}\cap B_{x_3})=n-2$ then
$\dim (B_{x_1}\cap B_{x_2})=n-2$ for a general $x_2$, a contradiction.
Thus $B_{x_1}$ and $B_{x_2}$ are either disjoint or they coincide for
  $(x_1,x_2)$ in an open subset of $X\times X$. 
The algebraic family of divisors $B_x$ thus determines a
morphism $p:X\to A$ to a smooth curve $A$. A general fiber of
$p$ is $B_x$.  Let  $Y$ be a general fiber of $p$
and  $x\in Y$  a general smooth point. 
We have an $(n-2)$-dimensional family of curves $C_t$
through $x$ and all of these are contained in $Y$.
$(C_t\cdot H)=1$ and the Picard number of $Y$ is 1 
as before. Thus $-K_Y\cong -K_X|_Y\cong nH|_Y$. 
Let $\pi:\bar Y\to Y$ denote the normalization of $Y$.
The curves $C_t$ lift to $\bar Y$,
and 
$$
-K_{\bar Y}=-K_Y+(\mbox{conductor of $\pi$}), 
$$
where the conductor of $\pi$ is an effective divisor
which is zero iff $\pi$ is an isomorphism (this is a special case
of duality for finite morphisms, 
cf.\ \cite[Ex.III.7.2]{hartsh}).
 Thus
$$
-(K_{\bar Y}\cdot \pi^*H^{n-2})\geq -(K_Y\cdot H^{n-2})
= n(H|_Y)^{n-1}=n(\pi^*H)^{n-1}. 
$$
So $\bar Y\cong \p^{n-1}$ by (\ref{p^nlem}). 
Furthermore, we also obtain that
$-K_{\bar Y}=n\pi^*H+(\mbox{conductor})$,
which implies that the conductor of $\pi$ is zero.
Thus 
$Y\cong \p^{n-1}$ and  $p:X\to A$
is generically a $\p^{n-1}$-bundle.
This implies that $p$ has a section (cf.\  \cite[X.6-7]{locfields}), hence
every fiber of $p$ has a reduced irreducible component.

Let $W\subset X$ be a reduced irreducible component of a fiber of $p$
and $x\in W$ a general smooth point. The above argument applies
also to $W$, and we obtain that $W\cong \p^{n-1}$ and 
$W$ is a connected component of its fiber.
This shows that $X$ is a $\p^{n-1}$-bundle over $A$. \qed
\end{say}

\begin{ack}  
Partial financial support was provided by  the NSF under grant numbers
 DMS-9800807 and 
DMS-9622394. 
\end{ack}

\vskip1cm

\noindent Johns Hopkins  University, Baltimore MD 21218

\begin{verbatim}kachi@math.jhu.edu\end{verbatim}

\noindent Princeton University, Princeton NJ 08544-1000

\begin{verbatim}kollar@math.princeton.edu\end{verbatim}

\end{document}